\documentclass[11pt]{article}

\usepackage{amsmath,amssymb,amsthm}
\usepackage{latexsym}
\usepackage{graphicx}
\usepackage{subfigure}
\usepackage{epsfig}

\newcommand{\Hmm}[1]{\leavevmode{\marginpar{\tiny%
$\hbox to 0mm{\hspace*{-0.2mm}$\leftarrow$\hss}%
\vcenter{\vrule depth 0.1mm height 0.1mm width \the\marginparwidth}%
\hbox to
0mm{\hss$\rightarrow$\hspace*{-0.5mm}}$\\\relax\raggedright #1}}}

\newcommand{\beeq}{\begin{equation}}
\newcommand{\eneq}{\end{equation}}

\renewcommand{\mod}{{\rm mod\ }}

\newcommand{\bear}{\begin{eqnarray}}
\newcommand{\eear}{\end{eqnarray}}
\newcommand{\beq}{\begin{equation}}
\newcommand{\eeq}{\end{equation}}

\newcommand{\spec}{{\rm spec}}

\newcommand{\R}{{\mathbb R}}

\newcommand{\Z}{{\mathbb Z}}

\newcommand{\calg}{\,{\mathfrak g}}

\newcommand{\calH}{{\mathcal H}}

\newcommand{\Laplace}{\triangle}

\def\calge1{\calg_{\vec{e_1}}}

\def\bm{\left[\begin{array}{cc}}
\def\endm{\end{array}\right]}
\def\ker{{\rm ker}}
\def\Ran{{\rm Ran}}
\def\Hil{{\mathcal H}}

\newtheorem{theorem}{Theorem}
\newtheorem{lemma}[theorem]{Lemma}

\newtheorem{claim}[theorem]{Claim}
\newtheorem{corollary}[theorem]{Corollary}
\renewcommand{\d}{\partial}
\newcommand{\<}{\langle}
\renewcommand{\>}{\rangle}

\theoremstyle{remark}

\def\norm[#1][#2]{\Vert #1 \Vert_{#2}}

\def\Util3{\tilde{U}^{(3)}}

\begin{document}

\title {Numerical verification of a gap condition for linearized NLS}
\author{L.\ Demanet, W.\ Schlag\thanks{The second author was partially
    supported by the NSF grant DMS-0300081 and a Sloan Fellowship. We
    would like to thank D. Pelinovsky and L. Ying for interesting
    discussions.}}

\maketitle

\begin{abstract}We make a detailed numerical study of the spectrum of two
  Schr\"{o}dinger operators $L_\pm$ arising in the linearization of
  the supercritical nonlinear Schr\"{o}dinger equation (NLS) about the
  standing wave, in three dimensions. This study was motivated by a
  recent result of the second author on conditional asymptotic
  stability of solitary waves in the case of a cubic nonlinearity.
  Underlying the validity of this result is a spectral condition on
  the operators $L_\pm$, namely that they have no eigenvalues nor
  resonances in the gap (a region of the positive real axis between
  zero and the continuous spectrum,) which we call the gap property.
  The present numerical study verifies this spectral condition, and
  further shows that the gap property holds for NLS
  exponents of the form $2 \beta + 1$, as long as $\beta_* < \beta \leq
  1$, where
\[
\beta_* = 0.913958905 \pm 1e-8.
\]

Our strategy consists of rewriting the original eigenvalue problem via
the Birman-Schwinger method. From a numerical analysis viewpoint, our
main contribution is an efficient quadrature rule for the kernel
$1/|x-y|$ in $\R^3$, i.e., provably spectrally accurate. As a result,
we are able to give similar accuracy estimates for all our eigenvalue
computations. We also propose an improvement of the Petviashvili's
iteration for the computation of standing wave profiles which
automatically chooses the radial solution.

All our numerical experiments are reproducible. The Matlab code can be
downloaded from
\begin{verbatim}http://www.acm.caltech.edu/~demanet/NLS/\end{verbatim}

\end{abstract}

\section{Introduction}

Suppose that $\psi(t,x)=e^{it\alpha^2}\phi(x)$ with $\alpha\ne0$
and $x\in\R^d$ is a standing wave solution of the NLS \beeq
i\partial_t\psi + \Delta\psi + |\psi|^{2\beta}\psi =0,
\label{eq:NLS} \eneq where $0<\beta<\frac{2}{d-2}$ if $d\ge3$ and
$0<\beta<\infty$ if $d=1,2$. Here we assume that
$\phi=\phi(\cdot,\alpha)$ is a ground state, i.e.,
\[
\alpha^2 \phi - \Delta\phi = \phi^{2\beta+1}, \qquad \phi>0.
\]
It is known that such $\phi$ exist and that they are radial, smooth,
and exponentially decaying, see Berestycki, Lions~\cite{BerLio}
and for uniqueness, Kwong~\cite{Kwo}. In one dimension $d=1$,
these ground states are explicitly given as
\beeq \label{eq:1dexp}
\phi(x) =
\frac{(\beta+1)^{\frac{1}{2\beta}}}{\cosh^{\frac{1}{\beta}}(\beta
x)} \eneq when $\alpha=1$ (for other values of $\alpha\ne 0$
rescale), but in higher dimensions no explicit expression is
known. From now on, we shall assume that $d=3$.

A much studied question is the stability of these standing waves,
both in the orbital (or Lyapunov) sense and the asymptotic sense.
For the former, see for example Grillakis, Shatah,
Strauss~\cite{GSS1}, \cite{GSS2}, Weinstein~\cite{Wei1},
\cite{Wei2}, Grillakis~\cite{Grill}, and for the latter, Buslaev,
Perelman~\cite{BP1}, Cuccagna~\cite{Cuc}. Reviews are in
Strauss~\cite{Str} and Sulem, Sulem~\cite{SulSul}.

In order to study stability, one generally linearizes around the
standing wave. This process leads to matrix Schr\"{o}dinger
operators of the form
$$
\calH = \calH_0 + V = \left[
\begin{array}{cc}
-\Delta + \alpha^2 & 0 \\
0 & \Delta - \alpha^2
\end{array}\right] +
\left[\begin{array}{cc} -V_1 & -V_2\\
V_2 & V_1
\end{array}
\right]
$$
on $L^2(\R^d)\times L^2(\R^d)$. Here,  $V_1 =
(\beta+1)\phi^{2\beta}$ and $V_2 = \beta \phi^{2\beta}$.

Conjugating $\Hil$ by the matrix $\bm 1 & i \\ 1 & -i\endm$ leads
to the matrix operator
\[ \bm 0 & i L_{-} \\ -iL_{+} & 0 \endm \]
with \begin{align*} L_- &= -\Laplace + \alpha^2 - \phi^{2\beta} \\
L_+ &= -\Laplace + \alpha^2 - (2\beta+1)\phi^{2\beta}
\end{align*}
The continuous spectrum of both $L_-$ and $L_+$ equals
$[\alpha^2,\infty)$. Since $L_-\phi=0$ and $\phi>0$, it follows
that zero is a simple eigenvalue and  the bottom of the spectrum
of $L_-$. Moreover, $L_+\partial_j \phi=0$ for $1\le j\le 3$ so
that $\ker(L_+)\subset \{\partial_j\phi\::\:1\le j\le 3\}$. In
fact, for monomial nonlinearities it is known that there is
equality here, see Weinstein~\cite{Wei2}\footnote{In this paper a
restriction $\beta\le1$ is imposed in $d=3$, but using Kwong's
results~\cite{Kwo} allows one to obtain the full range $\beta<2$
by means of Weinstein's arguments}, and that there is a unique
negative bound state of $L_+$.

It is known that  $L_-\ge0$ implies that the spectrum $\spec(\calH)$
satisfies $\spec(\calH) \subset \R \cup i \R$ and that all points of
the discrete spectrum other than zero are eigenvalues whose
geometric and algebraic multiplicities coincide. On the other hand,
the zero eigenvalue of $\Hil$ has geometric multiplicity four and
algebraic multiplicity eight provided $\beta\ne\frac{2}{3}$, whereas
for the $L^2$-critical case $\beta=\frac23$ the algebraic
multiplicity increases to ten. For this see~\cite{Wei1}, \cite{BP1}
or~\cite{RSS1}, \cite{ErdSch}.

In order to carry out a meaningful asymptotic stability analysis
it is essential to understand the discrete spectrum of $\Hil$. The
root space at zero was completely described by
Weinstein~\cite{Wei2}. Moreover, it is also well-known that
\[ \spec(\Hil) \subset \R \text{\ \ iff\ \ }\beta\le\frac23 \]
whereas in the range $\frac23<\beta <2$ there is a unique pair of
simple complex-conjugate eigenvalues $\pm i\gamma$. This latter
property reflects itself in the nonlinear theory in the following
way: Orbital stability holds iff $\beta<\frac23$, see Berestycki,
 Cazenave~\cite{BerCaz},  Weinstein~\cite{Wei2}, Cazenave,
Lions~\cite{CazLio}, and~\cite{GSS1}, \cite{GSS2}.

In \cite{schlag} the second author investigated conditional
asymptotic stability for the unstable case $\beta=1$. This
analysis depended on the fact that zero is the only eigenvalue of
$\Hil$ in the interval $[-\alpha^2,\alpha^2]$ and that the edges
$\pm\alpha^2$ are not resonances. The fact that $\pm \alpha^2$ are
neither eigenvalues nor resonances is the same as requiring that
the resolvent $(\Hil - z)^{-1}$ remains bounded on suitable
weighted $L^2(\R^3)$ spaces for $z$ close to $\pm\alpha^2$.

Using some ideas of Perelman~\cite{Pe2}, it is shown
in~\cite{schlag} that these properties can be deduced from the
following properties of $L_+, L_-$: {\bf Neither $L_+$ nor $L_-$
have any eigenvalues in the gap $(0,\alpha^2]$ and $L_-$ has no
resonance at $\alpha^2$.}

In one dimension $d=1$, the spectral properties of $L_-$ and $L_+$
can be determined completely since the generalized eigenfunctions of
these operators (more precisely, the Jost solutions) can be given
explicitly in terms of certain hypergeometric functions, see
Fl\"ugge~\cite{Flug}, Problem~39 on page~94. This is due to the
special form of the ground state~\eqref{eq:1dexp}.

Unfortunately, it seems impossible to determine similar properties
for the case of three dimensions by means of purely analytical
methods. We therefore verify this gap property of $L_{\pm}$
numerically via the Birman-Schwinger method. We will refer to the
{\bf gap property} as the fact that {\bf $L_{\pm}$ have no
eigenvalues in $(0,\alpha^2]$ and no resonance at $\alpha^2$.} Our
main result is as follows.

\begin{claim}
There exists a number $\beta_*=.913958905 \pm 1e-8$  so that for all
$\beta_*<\beta\le1$ the gap property holds.
\end{claim}

This statement can also be continued beyond $\beta=1$. We only went
up to $1$ since $\beta=1$ alone is needed in~\cite{schlag}. In the
range $\beta < \beta_*$, our numerical analysis shows that the
operator $L_+$ has eigenvalues in the gap $(0,1]$. This is perhaps
surprising, since it shows that the gap property does not hold for
the entire $L^2(\R^3)$ super-critical range $\frac23<\beta<2$.
However, it does hold at $\beta=1$. In particular, the method of
proof from~\cite{schlag} does not apply to all $\beta>\frac23$ since
it relies on the gap property.  In contrast, in one dimension $d=1$,
Krieger and the second author showed that this method does apply to
the entire super-critical range $\beta>2$. In fact, there, the gap
property does hold for all $\beta>1$, see~\cite{Flug}.

In the remainder of this paper, we will be concerned with the
description of the numerical method, and will study its convergence
properties. Note that we did not formulate our main result as a
theorem because the error bound we give on $\beta_*$ is based
on numerical observations of convergence of the method. A fully
rigorous justification would involve either an \emph{a priori} estimate
for this bound, or a provably reliable \emph{a posteriori}
estimate. For this, we would probably need to (1) derive quantitative
estimates of regularity and decay for $\phi(x)$ and other functions,
with explicit values of the constants, and (2) give a full treatment of the
propagation of round-off errors from a countless number of sources.
This type of heroic exercise would contribute little to the
understanding of the accuracy of the numerical method, so we chose to
spare the reader.

\section{Description of the numerical method}\label{sec:num}

\subsection{The Birman-Schwinger method}\label{sec:BS}

A direct numerical computation of the eigenvalues of $L_\pm$ would be
problematic near $\alpha^2$, the edge of the continuous spectrum.  For
example, it is unclear if a perceived numerical eigenvalue at $0.99
\, \alpha^2$ belongs to the gap or if it has escaped from the continuous
spectrum. Decay of the corresponding eigenfunction might help in
making a decision, but this criterion is unacceptable over a
truncated computational domain.  Instead, we will reformulate the
problem by the Birman-Schwinger method, which we now recall.

Let $H=-\Laplace - V$, where $V>0$ is a bounded potential that
decays at infinity. In our case, this is $H = L_\pm - \alpha^2 I$.
We would like to filter out positive eigenvalues, so assume
$Hf=-\lambda^2 f$ where $\lambda>0$ and $f\in L^2$. Then $g=Uf$,
where $U = \sqrt{V}$ satisfies
\[ g = U (-\Laplace + \lambda^2)^{-1} U g.  \]
In other words, $g\in L^2 $ is an eigenfunction of
\[ K(\lambda) = U (-\Laplace + \lambda^2)^{-1} U \]
with eigenvalue one. Note that $K(\lambda)$ is a compact, positive
operator. Conversely, if $g\in L^2$ satisfies $K(\lambda)g=g$,
then
\[ f:=U^{-1}g = (-\Laplace + \lambda^2)^{-1} U g \in L^2\] and
$Hf=-\lambda^2 f$. Moreover, the eigenvalues of $K(\lambda)$ are
strictly increasing as $\lambda\to 0$. Hence, we conclude that
\[ \# \Big\{ \lambda\::\: \ker(H-\lambda^2)\ne \{0\} \Big\}
=\# \Big\{ E>1\::\: \ker(K(0)-E)\ne \{0\} \Big\} \; ,
\]
counted with multiplicity.

Finally, in view of the symmetric resolvent identity, viz.
\[ (H-z)^{-1} = (-\Laplace-z)^{-1} +
(-\Laplace-z)^{-1}U\Big[I-U(-\Laplace-z)^{-1}U\Big]^{-1}
U(-\Laplace-z)^{-1}.
\]
This shows that the Laurent expansion of $(H-z)^{-1}$ around $z=0$
does not involve negative powers of $z$ iff $I+U(-\Laplace-z)^{-1}U$
is invertible at $z=0$ which is the same as requiring that
\[\ker\{I-U(-\Laplace)^{-1}U\}=\{0\}\] because of the Fredholm
alternative (assuming that $V$ decays sufficiently fast at infinity to
insure compactness). In other words, if $H$ has no resonance or
eigenvalue at the origin, then $K(0)$ will not show an eigenvalue $E =
1$, and conversely.

Let us count the eigenvalues $\{\lambda_j\}_{j=1}^\infty$ of
$K(0)=U(-\Laplace)^{-1}U$ (which are all non-negative) in decreasing
order. Then we arrive at the following conclusion: Let $N$ be a
positive integer.  {\bf Then the  operator $H$ has exactly $N$
negative eigenvalues and neither an eigenvalue nor a resonance at
zero iff $\lambda_1\ge \ldots \ge\lambda_N>1$ and
$\lambda_{N+1}<1$.}

Note that spectrum of the self-adjoint, compact Birman-Schwinger
operator $K(0)$ is discrete and robust to numerical perturbations
near $E = 1$, which is the desired numerical effect. Since $K(0)$ is
compact, eigenvalues cannot accumulate to $E = 1$ from below, and
are in no way related to the continuous spectrum of~$H$.

\bigskip In view of the preceding, we therefore need to show the
following to justify~\cite{schlag}: For $\beta=1$, the second
largest eigenvalue of\footnote{Note that due to the scaling
$x\mapsto \alpha x$ and
$\phi(x,\alpha)=\alpha^{\frac{1}{\beta}}\phi(\alpha x)$ we may
assume that $\alpha=1$. From now on we set $\phi(x)=\phi(x,1)$.}
\[ K_{-}(x,y) = \frac{\phi^{\beta}(x)\phi^\beta(y)}{4\pi|x-y|} \]
is below one, and the fifth largest eigenvalue of
\[ K_{+}(x,y) = (2\beta+1)\frac{\phi^{\beta}(x)\phi^\beta(y)}{4\pi|x-y|} \]
is below one. These properties will then imply the gap property,
i.e., $L_{\pm}$ have no eigenvalues in $(0,1]$ and no resonance
at $1$.

\subsection{The modified Petviashvili's iteration}\label{sec:petv}

The first step of the numerical method is to find the soliton
$\phi(x)$, which is the unique positive, radial, decaying solution
of
\begin{equation}\label{eq:NLell}
- \Delta \phi + \phi = \phi^{2 \beta+1},
\end{equation}
unique up to translation. As mentioned earlier, $\phi(x)$ is in fact
exponentially decaying. A naive approach would be to solve a descent
equation like
\[
\frac{\d u}{\d t} = \Delta u - u + |u|^{2 \beta}u,
\]
but, as shown in \cite{BerLioPel}, this equation is unstable near the
fixed manifold of interest. Instead, we will solve the
\emph{modified Petviashvili's iteration} which reads
\begin{equation}\label{eq:petv}
\phi_{n+1} = M_n^\gamma (I - \Delta)^{-1}(|\phi_n|^{2 \beta} \phi_n) + \delta \sum_{j = 1}^3 R_{n,j} \frac{\d \phi_n}{\d x_j}
\end{equation}
The initial guess $\phi_0$ can for example be taken as a Gaussian. The
choice of constants $M_n, R_{n,j}, \gamma$ and $\delta$ is crucial for
convergence of the iteration, and is given by
\begin{align}\label{eq:petvconst}
M_n &= \frac{\int (1 + |\xi|^2) (\widehat{\phi_n})^2 \, d\xi}{\int \widehat{\phi_n}\,\widehat{(|\phi_n|^{2 \beta}\phi_n)} \, d\xi}, \\
R_{n,j} &= \frac{\int (1 + |\xi|^2) \widehat{\phi_n} \,\widehat{\d_j\phi_n} \, d\xi}{\int \widehat{\d_j \phi_n}\,[\d_j (|\phi_n|^{2 \beta}\phi_n)]^{\wedge} \, d\xi}, \\
\gamma &= \frac{2 \beta + 1}{2 \beta}, \\
\delta &= -1/2,
\end{align}
where $\d_j = \frac{\d}{\d x_j}$, and the hat denoting Fourier
transformations.  This iteration, without the second term, was
introduced by Petviashvili in 1976, and convergence was proved
recently in \cite{PelSte}. The addition of the second term is a minor
increment whose purpose is to fix a potential source of instability
due to numerical discretization and to force the iteration to choose
the radial soliton (centered at the origin).  This will be explained
and justified in section \ref{sec:conv}.

Numerically, $(I-\Delta)^{-1}$ is realized in the Fourier domain, and
Fourier transformations are implemented via the Fast Fourier Transform
(FFT).  Discretization issues are addressed in the next section. In
practice, the iteration is accelerated via Aitken's method applied
pointwise, i.e.,
\[
\phi^A_{n}(x) = \phi_n(x) - \frac{(\phi_{n+1}(x) - \phi_n(x))^2}{\phi_{n+2}(x) - 2\phi_{n+1}(x)+\phi_n(x)}.
\]

The iteration is stopped when the Euler-Lagrange equation
(\ref{eq:NLell}) is satisfied up to some very small tolerance $\tau$
in $L^2$.  The resulting approximation of the soliton will be denoted by
$\tilde{\phi}$.

\subsection{Truncation and Discretization}

We now take up the task of computing the eigenvalues of the
Birman-Schwinger operators as defined in Section~\ref{sec:BS}. The
first step is to truncate the three-dimensional computational domain
to a cube of sidelength $L$ centered at the origin and to
discretize functions $f(x)$ by evaluating them on the regular grid
\begin{equation}\label{eq:grid}
x_j = (j_1, j_2, j_3) \frac{L}{N},
\end{equation}
with $j_1, j_2, j_3$ integers obeying $-N/2 \leq j_k \leq N/2 - 1$.
Operators are, in turn, discretized as matrices acting on `vectors' of
function samples $f(x_j)$. Tools of numerical linear algebra can then
be invoked to compute the eigenvalues of these matrices.

Typical values of $L$ and $N$ for which discretizing $K_\pm$ is
expected to be reasonably accurate are $L \simeq 20$ and $N \simeq
100$.  In this context, several vectors of $N^3 \simeq 10^6$ function
samples can comfortably be stored simultaneously in the memory of a
2005-era computer, but we cannot yet afford to manipulate matrices
containing $N^6 \simeq 10^{12}$ elements. This rules out the
possibility of using popular approaches such as the QR algorithm,
which compute eigenvalues by operating directly on the matrix entries.

Instead, we will resort to a modification of the power method, known
as the implicitly restarted Arnoldi iteration, which is
implemented in Matlab's \emph{eigs} command \cite{eigs}. This method has the
advantage of only requiring applications of the operator to
diagonalize, i.e., matrix-vector products. For well-conditioned
problems, such as the one we are addressing, \emph{eigs} computes the
top eigenvalues of the finite matrix up to machine precision, i.e.,
about 15 decimal digits in Matlab.

The only remaining issue is then to find a good discretization
$\tilde{K}_\pm$ of the Birman-Schwinger operators $K_\pm$, and to
quantify the accuracy. Multiplication by $\phi(x)$ to some power
will be done sample-wise on the grid $x_j$. Inverting minus the
Laplacian in a space of decaying functions over $\mathbb{R}^3$, or
equivalently convolving with the fundamental solution $G(x) =
\frac{1}{4 \pi |x|}$, is a bit more complicated. Discretizing $G(x)$
by sampling at $x_j$ is quite inaccurate and is problematic
if $x = 0$ belongs to the grid $x_j$. Dividing by $|\xi|^2$ in
frequency poses similar difficulties. Instead, for reasons which
will be explained in Section~\ref{sec:conv}, we use the following
discretization,
\begin{equation}\label{eq:discfun}
\tilde{G}(x_j) = \left\{
\begin{aligned}
  \frac{1}{2 \pi^2 |x_j|} \left( \frac{L}{N} \right)^3 \, \text{Si}(\frac{\pi N |x_j|}{L}) \qquad &\text{if } x_j \ne 0, \\
  \frac{1}{2 \pi} \left( \frac{L}{N} \right)^2 \qquad &\text{if } x_j = 0,
\end{aligned}
\right.
\end{equation}
where $\text{Si}(x) = \int_0^x \frac{\sin t}{t} \, dt$. The discrete
(circular) convolution of $\tilde{G}(x_j)$ with a vector $f(x_j)$ is
computed efficiently as the multiplication $\widehat{\tilde{G}(x_j)}
\, \widehat{f(x_j)}$ of their respective FFT, followed by an inverse FFT.

In this context, applying the full Birman-Schwinger operator, say $K_-$, to a
vector of samples $f(x_j)$ consists of the obvious sequence of steps:
(1) multiply $f(x_j)$ by $\tilde{U}(x_j) = \tilde{\phi}(x_j)^\beta$, (2) perform an
FFT, (3) multiply the result by $\widehat{\tilde{G}(x_j)}$, (4)
perform an inverse FFT, and finally (5) multiply the result by
$\tilde{U}(x_j)$ again.

Since the complexity of a one-dimensional FFT in Matlab is $O(N \log
N)$ operations for most values of $N$ (not necessarily a power of
two), one application of $\tilde{K}_\pm$ will require $O(N^3 \log N)$
operations. This is a substantial improvement over the naive
matrix-vector product which would require $O(N^6)$ operations.

\section{Convergence analysis}\label{sec:conv}
\subsection{The modified Petviashvili's iteration}\label{sec:petv_conv}

In this section we discuss convergence of the iteration
(\ref{eq:petv}). The first result in this direction, in the case
$\delta = 0$, can be found in ref.~\cite{PelSte}. To make this
discussion self-contained, we recall their argument and apply it to
our specific problem.

\begin{theorem}\label{teo:petv}{\bf (Pelinovsky, Stepanyants.)}
  Let $\phi(x)$ be the unique radial solution of (\ref{eq:NLell}) and
  $H^1_{r}(\R^3)$ denote the subset of all radial functions in
  $H^1(\R^3)$. Consider the iteration (\ref{eq:petv}) with $M_n$ and
  $\gamma$ given by equations (\ref{eq:petvconst}), but $\delta = 0$.
  Then there exists an open neighborhood $\mathcal{N}$ of $\phi$ in
  $H^1_r(\R^3)$, in which $\phi$ is the unique fixed point and
  (\ref{eq:petv}) converges to $\phi$. The iteration is strictly
  stable in the sense that, for all $\phi_0 \in \mathcal{N}$,
\begin{equation}\label{eq:strictstab}
|| \phi_{n+1} - \phi ||_1 \leq (1-C) || \phi_{n} - \phi ||_1, \qquad 0 < C \leq 1,
\end{equation}
where $|| \cdot ||_1$ is the norm in the Sobolev space $H^1(\R^3)$.
\end{theorem}
\begin{proof}

  Put $p = 2 \beta +1$. Let us first write down the linearized
  iteration, about the fixed point $\phi$, for the perturbation $w_n =
  \phi_n - \phi$. It reads\footnote{Throughout this paper, we use the following convention for the Fourier transform:
\[
\hat{f}(\xi) = \int e^{-ix\cdot\xi} f(x)\, dx, \qquad f(x) =
\frac{1}{(2 \pi)^3} \int e^{i x\cdot\xi} \hat{f}(\xi)\, d\xi.
\]}
\begin{equation}\label{eq:linpetv}
\hat{w}_{n+1} = (1-p) \gamma a_n \hat{\phi} + p \frac{\widehat{\phi^{p-1}} *
\hat{w}_{n}}{1 + |\xi|^2},
\end{equation}
where $*$ denotes convolution, $a_n$ comes from the
linearization of $M_n$, and is given by
\[
a_n = \frac{\int w_n \phi^p}{\int \phi^{p+1}}.
\]
(Up to higher-order terms, we have the asymptotic relation $1+(1-p)a_n
\sim M_n$.) This formula suggests that $\phi$ plays a special role for
the stability of the linearized iteration. Indeed, we claim that we
can actually expand $w_n$ as
\begin{equation}\label{eq:expwn}
w_n = a_n \phi + q_n,
\end{equation}
where $a_n$ is exactly as defined above, and $q_n$ is a remainder. In
order to see this, let us introduce the operator $A = (I -
\Delta)^{-1} \mathcal{H}$, where $\mathcal{H} = I - \Delta - p
\phi^{p-1}$. It is easy to check that is it bounded and self-adjoint
with respect to the $H^1$ inner product,
\[
(f,g) \equiv \< f, (I-\Delta) g\>.
\]
The operator $A$ therefore provides a spectral decomposition of
$L^2(\R)$, orthogonal with respect to the inner product
$(\cdot,\cdot)$. It was noticed in \cite{PelSte}, (or by a
straightforward extension of their argument,) that the spectra of $A$
and $\mathcal{H}$ obey
\begin{align*}
\mbox{dim}(\mbox{neg}(A)) &= \mbox{dim}(\mbox{neg}(\mathcal{H})) = 1, \\
\mbox{null}(A) &= \mbox{null}(\mathcal{H}) = 3.
\end{align*}
The first four eigenfunctions of $A$ are precisely $\phi$ and $\d_i
\phi$, with eigenvalues $1-p < 0$ and $0$ respectively. Equation
(\ref{eq:expwn}) is just the expansion of $w_n$ in this orthogonal
system. Since the iterates $\phi_n$ are all radial, the components
along $\d_i \phi$ are zero. The remainder $q_n$ belongs to the space
$Y_p$ defined by
\[
Y_p = \{ u \in L^2(\R^3): \< u, \phi^p \> = \< u, \d_i \phi^p \> = 0 \}.
\]
In the space $Y_p$, the spectrum of $A$ is strictly positive,
bounded from below by the fifth eigenvalue\footnote{If there is no
fifth eigenvalue, then $\lambda_5$ equals the edge of the essential
spectrum} $\lambda_5 > 0$ and from above by
\[
\lambda_M = \sup_u \frac{(u, A u)}{(u,u)} = 1 - p \, \inf_u
\frac{\<u, \phi^{p-1} u \> }{(u,u)} = 1.
\]
For the last equality we have used the fact that $\phi(x) > 0$ and
$\phi(x)\to0$ as $|x|\to\infty$.

The recurrence equations for $a_n$ and $q_n$ can be found
from equation (\ref{eq:linpetv}),
\begin{align*}
a_{n+1} &= (p - \gamma (p - 1)) a_n, \\
q_{n+1} &= (I - A) q_n.
\end{align*}
The choice we made for $\gamma$ ensures that the component along
$\phi$ is immediately put to zero (in the linearized iteration.) It is
also clear that we have $|| q_{n+1} ||_1 \leq (1 - \lambda_5) || q_n
||_1$ in $H^1(\R^3)$. In the scope of the linearized iteration,
equation (\ref{eq:strictstab}) follows with $C = \lambda_5$.

Call $\mathcal{P}$ the nonlinear operator for the Petviashvili
iteration in the case $\delta = 0$, so that (\ref{eq:petv}) is written
as $\phi_{n+1} = \mathcal{P}(\phi_n)$. We follow \cite{PelSte} and
apply the contraction mapping theorem in a neighborhood of the pixed
point $\phi$, within the closed subspace $H^1_r(\R^3)$. We have
already computed the Fr\'{e}chet derivative of $\mathcal{P}$ at
$\phi$,
\begin{equation}\label{eq:frechetP}
\mathcal{P}'(\phi)w = (I - A)P_{Y_p}w,
\end{equation}
where $P_{Y_p}$ is the projection onto $Y_p$, orthogonal in $H^1$. Let
$0 < \epsilon < \lambda_5$. By continuity of $\mathcal{P}'(u)$ as a
function of $u$, in the operator $H^1$ norm, we can assert that there
exists a small open neighborhood $\mathcal{N}$ of $\phi$ in which
\[
|| \mathcal{P}'(u) ||_{1 \to 1} < 1 - \lambda_5 + \epsilon.
\]
By a standard application of the contraction mapping theorem (see
\cite{HutPym} p.126), the fixed point is unique in $\mathcal{N}$ and
we have the estimate
\[
|| \phi_{n+1} - \phi ||_1 \leq (1 - \lambda_5 + \epsilon) || \phi_n -
\phi||_1.
\]
This concludes the proof.
\end{proof}

Let us now examine how the above argument generalizes to the case
$\delta \ne 0$. The purpose of the second term in equation
(\ref{eq:petv}) is precisely to put to zero the components along the
three basis functions $\d_i \phi$, should the initial condition not be
radial. This is also useful in the context of the numerical
realization of $\mathcal{P}$, since numerical round-off errors do not
correspond in general to radial perturbations. To be precise, the
linearized iteration (\ref{eq:linpetv}) becomes
\[
\hat{w}_{n+1} = (1-p) \gamma a_n \hat{\phi} + 2 \delta \sum_{j =
1,2,3} b_{n,j} \widehat{\d_j \phi} + p \frac{\widehat{\phi^{p-1}} *
\hat{w}_{n}}{1 + |\xi|^2},
\]
with $a_n$ as previously and
\[
b_{n,i} = \frac{\int w_n \d_i \phi^p}{\int \d_i \phi \d_i \phi^p}
\]
is the component of $w_n$ along $\d_i \phi$, in the natural $H^1$ inner
product. So we have
\[
w_n = a_n \phi + \sum_j b_{n,j} \d_j \phi + q_n,
\]
and the recurrence relation for $b_{n,i}$ is
\[
b_{n+1,i} = (1 + 2 \delta) b_{n,i}.
\]
Choosing $\delta = -1/2$ as advocated previously will put to zero the
components along $\d_i \phi$.

In practice the modification $\delta = -1/2$ works very well, see
section \ref{sec:results}. The analysis of the linearized iteration
does not pose any difficulty. However, we have been unable to extend
the argument of theorem \ref{teo:petv} to the full \emph{nonlinear}
iteration. The operator $\mathcal{P}$ is in general not defined on
$H^1$, because the constants $R_{n,i}$ involve 3/2 derivatives of
$\phi_n$ in $L^2$. A fortiori, $\mathcal{P}$ is not
Frechet-differentiable for functions in $H^1$. Note that the
contraction argument cannot work, since the soliton $\phi$ is not
unique in $H^1(\R^3)$ -- it is unique up to a translation. The question
of stability of the modified Petviashvili's iteration in the
non-radial case is possibly related to the problem of proving
uniqueness of the radial soliton, which, in itself, is not trivial.

\subsection{Truncation and Discretization}

In this section we prove \emph{spectral convergence} of the proposed
discretization of the Birman-Schwinger operator, which is numerical
analysis jargon for (almost) exponential convergence with respect to
the large discretization parameters $L$ and $N/L$. Recall that we
denote by $x_j = (j_1,j_2,j_3)\frac{L}{N}$, $j \in \Z^3$, the nodes of
a cubic grid with spacing $L/N$ in all three directions, not
necessarily bounded in some of the arguments that follow. Obviously,
we take $L \leq N$. Unless otherwise specified, we are considering the
operator $K_-$, since $K_+ = (2 \beta + 1)K_-$.

We now formulate our main result. Most of the rest of this section is
devoted to its justification. We will use the notation $\< x \> =
\sqrt{1+|x|^2}$.

\begin{theorem}\label{teo:main}
  Let $K$ be the Birman-Schwinger operator as above, and $\tilde{K}$
  its numerical realization, extended to functions of continuous $x$
  by sampling and interpolation (see below for details). Let $\delta >
  0$ be the decay rate of $U(x) = \phi(x)^\beta$, $|U(x)| \leq C \cdot
  e^{-\delta |x|}$.  Assume that Petviashvili's method gives an accurate
  approximation $\tilde{\phi}$ of the soliton in the sense that, for
  some $\epsilon > 0$, and denoting $\tilde{U} = \tilde{\phi}^\beta$,
\begin{equation}\label{eq:assum_petv}
|U(x_j) - \tilde{U}(x_j)| \leq C \cdot \min(\frac{\epsilon}{\< x_j \>}, e^{-\delta |x_j|}).
\end{equation}
Then we have, for all $s > 3/2$ and $f \in H^s(\R^3)$, in exact arithmetic,
\begin{equation}\label{eq:main}
|| (K - \tilde{K}) f ||_{L^2} \leq C_{s} \cdot \left[ \epsilon + L e^{-\delta L/4} + \left( \frac{N}{L} \right)^{-s} \right] \cdot || f ||_{H^s},
\end{equation}
for some constant $C_s > 0$ depending on $s$.
\end{theorem}

Discretization is error-free in the context of the Shannon sampling
theorem. Let us introduce $B_{N/L}(\R^3)$, the space of band-limited
square-integrable functions,
\[
B_{N/L}(\R^3) = \{ u \in L^2(\R^3): \hat{u}(\xi) = 0, \; | \xi | > \pi N/L \}.
\]
The hat denotes Fourier transformation. We can then define the
\emph{sampling} operator $S$, and \emph{interpolation} operator $T$,
as
\[
S: B_{N/L} \to \ell^2, \quad f(x) \mapsto \{ f(x_j) \}.
\]
\[
T: \ell^2 \to B_{N/L}, \quad \{ f(x_j) \} \mapsto \sum_{j \in \Z^3} h(x-x_j) f(x_j).
\]
Here $h(x)$ is the interpolating kernel defined by $\hat{h}(\xi) =
\left( \frac{L}{N} \right)^3$ if $|\xi| \leq \frac{\pi N}{L}$, and
zero otherwise. The content of Shannon's sampling theory is that $S$
is an isometry from $B_{N/L}$ to $\ell^2$, and $T$ is in that
context both the adjoint and a left inverse for $S$ on its range
(hence the interpolation property of $h$.) Note that the properly
normalized $\ell^2$ norm is
\[
\| \{ f(x_j) \}_j \|_{\ell^2} = \sqrt{ \left( \frac{L}{N} \right)^3
\sum_j |f(x_j)|^2}
\]

With a slight abuse of notations, let us denote by $\tilde{G}$ the
operator of discrete convolution by $\tilde{G}(x_j)$, defined in
equation (\ref{eq:discfun}). This particular expression for the
weights $\tilde{G}(x_j)$ is chosen so that the inversion of the
Laplacian is \emph{exact} on decaying band-limited functions. By
`decaying', we mean a member of some weighted $L^2$ space. More
generally, let us introduce weighted Sobolev spaces as
\[
H^s_m(\R^3) = \{ u: \< x \>^m u(x) \in H^s(\R^3) \},
\]
and equipped with the norm
\[
|| u ||_{s,m} = || \< x \>^m u ||_{s}.
\]
Of course, $L^2_m = H^0_m$.

\begin{lemma}\label{teo:exact_inverse_laplacian}
We have
\[
(T \, \tilde{G} \, Sf)(x) = \frac{1}{4 \pi} \int_{\R^3}
\frac{1}{|x-y|} f(y) \, dy,
\]
for all $f \in B_{N/L}(\R^3) \cap L^2_{1}(\R^3)$.
\end{lemma}
\begin{proof}
  The restriction $f \in L^2_{1}$ is simply chosen so that the
  integral makes sense. Since $f \in B_{N/L}$, we can express it as
\[
f(x) = \sum_j h(x-x_j) f(x_j).
\]
Substituting in the integral, we get
\begin{align*}
(- \Delta)^{-1} f(x_j) &= \sum_{k \in \Z^3} \frac{1}{4 \pi} \int \frac{1}{|x_j-y|}  h(y-x_k) \, dy \; f(x_k) \\
&= \sum_{k \in \Z^3} w(j,k) f(x_k).
\end{align*}
The quadrature weights $w(j,k)$ can be computed explicitly. By Parseval,
\[
w(j,k) = \frac{1}{(2 \pi)^3} \int \frac{1}{| \xi |^2} e^{-i
  (x_j-x_k) \cdot \xi} \hat{h}(\xi) \, d\xi.
\]
In the event $x_j \ne x_k$, we can express this integral in a
spherical coordinate frame whose z-axis is aligned with $x_j - x_k$.
It becomes
\begin{align*}
  w(j,k) &= \frac{1}{(2 \pi)^2} \left( \frac{L}{N} \right)^3 \int_0^{\pi N/L} r^2 dr \, \int_0^\pi \sin \theta d\theta \, \frac{1}{r^2} e^{-i | x_j - x_k | r \cos \theta}, \\
  &= \left( \frac{L}{N} \right)^3 \frac{1}{(2 \pi)^2} \int_0^{\pi N/L} dr \, \int_{-1}^1 e^{i | x_j - x_k | r x} \, dx, \\
  &= \left( \frac{L}{N} \right)^3 \frac{1}{(2 \pi)^2}\, 2 \int_0^{\pi N/L} \frac{\sin | x_j - x_k | r}{| x_j - x_k | r}\, dr, \\
  &= \left( \frac{L}{N} \right)^3 \frac{1}{2 \pi^2} \frac{1}{ | x_j - x_k |}
  \mbox{Si}\left( \frac{\pi N | x_j - x_k |}{L} \right).
\end{align*}
Si$(x)$ is the special function $\int_0^x \sin t / t \, dt$. Thus we
have $w(j,k) = \tilde{G}(x_j - x_k)$, as in equation
(\ref{eq:discfun}). The special case $x_j = x_k$ can be handled
separately in a similar fashion, and one can check that it corresponds
to the limit $x_j - x_k \to 0$ in the general expression.

We have just proved that $\tilde{G} S f = S (-\Delta^{-1}) f.$ Now
the output $(- \Delta)^{-1} f$ is not in general square-integrable,
see also Lemma~\ref{teo:JK} below, but it is still band-limited. The
interpolation operator $T$ is still the inverse of sampling for
band-limited functions with some growth, so we obtain the desired
conclusion.
\end{proof}

{\bf Remark.} It is interesting to notice that the choice we made
for $\hat{h}$, the indicator of a ball, is not standard. The most
natural frequency window associated to a Cartesian grid is the
indicator of the cube $[-\pi N/L, \pi N/L]^3$. This choice of window
for the definition of $B^{\mbox{std}}_{N/L}$, the space of $L^2$
band-limited functions, would have made the sampling operator $S$
not only an isometry, but a \emph{unitary} map from
$B^{\mbox{std}}_{N/L}$ \emph{onto} $\ell^2$. In our context, with
the spherical window, $S$ is only unitary from $B_{N/L}$ to a
\emph{subset} of $\ell^2$, namely the range $S B_{N/L}$. However, we
do not believe the cubic window would have allowed us to formulate a
closed-form expression for the weights $\tilde{G}(x_j)$.

\medskip The functions we have to discretize, like the soliton
$\phi(x)$, are unfortunately not band-limited. The smoother the
function, the more accurate its sampling, however. The following
elementary lemma establishes this property.

\begin{lemma}\label{teo:sampling}
For all $s > 3/2$ and $f \in H^s(\R^3)$,
\[
\| f - TSf \|_{L^2} \leq C_s \left( \frac{N}{L} \right)^{-s} \, \| f \|_{s}.
\]
The constant $C_s$ is a decreasing function of $s \in (3/2, \infty)$.
\end{lemma}
\begin{proof}
  Sampling in $x$ corresponds to periodizing in $\xi$, so the Fourier
  transform of $TSf - f$ is
\[
F(\xi) = \sum_{m \in \Z^3} \hat{f}(\xi - 2 \pi m \frac{N}{L})
\chi_{|\xi| \leq \pi N/L} (\xi) - \hat{f}(\xi).
\]
We get two contributions, call them $(I)$ and $(II)$, in
\begin{align*}
\| F \|_2^2 = &\int_{|\xi|\leq \pi N/L} |\sum_{m \ne 0} \hat{f}(\xi - 2 \pi m \frac{N}{L})|^2 \, d\xi \\
&+ \int_{|\xi|> \pi N/L} |\hat{f}(\xi)|^2 \, d\xi.
\end{align*}
The second integral $(II)$ is easily bounded by
\[
(II) \leq \sup_{|\xi| > \pi N/L} \< \xi \>^{-2s} \, \cdot \, \| f \|_s^2
\leq \pi^{-2s} \left( \frac{N}{L} \right)^{-2s} \, \| f \|_s^2.
\]
The first integral $(I)$ can be expressed as
\[
(I) = \sum_{m \ne 0} \sum_{m' \ne 0} I_{m,m'}
\]
where
\begin{align*}
  I_{m,m'} &= \int_{B_0(\pi N/L)} |\hat{f}(\xi - 2 \pi m \frac{N}{L}) \hat{f}(\xi - 2 \pi m' \frac{N}{L})|, \\
&\leq \sup_{|\xi|\leq \pi N/L} \< \xi - 2 \pi m \frac{N}{L} \>^{-s} \, \cdot \, \sup_{|\xi|\leq \pi N/L} \< \xi - 2 \pi m' \frac{N}{L} \>^{-s} \\
&\qquad \times \int_{B_0(\pi N/L)} |\hat{f}(\xi - 2 \pi m \frac{N}{L})| \< \xi - 2 \pi m \frac{N}{L} \>^s \cdot |\hat{f}(\xi - 2 \pi m' \frac{N}{L})| \< \xi - 2 \pi m' \frac{N}{L} \>^s \, d\xi.
\end{align*}
Obviously $\sup_{|\xi|\leq \pi N/L} \< \xi - 2 \pi m \frac{N}{L} \>^{-s} \leq C \cdot \left( |m| \frac{N}{L} \right)^{-s}$, and the integral can be bounded by Cauchy-Schwarz. We have thus separated
\[
I_{m,m'} \leq C \cdot J_m J_{m'} \cdot \left(\frac{N}{L} \right)^{-2s},
\]
where
\[
J_m = |m|^{-s} \sqrt{\int_{|\xi + 2 \pi m \frac{N}{L}|\leq \pi N/L} |\hat{f}(\xi)|^2 \< \xi \>^{2s} \, d\xi}.
\]
Each point $\xi$ in the frequency domain is covered by at most one
ball $B_{2 \pi m \frac{N}{L}}(\pi N/L)$. Consequently, Cauchy-Schwarz
for sequences gives
\[
\sum_{m \ne 0} J_m \leq C \cdot \sqrt{\sum_{m \ne 0} |m|^{-2s}} \cdot
\sqrt{\int_{\R^3} |\hat{f}(\xi)|^2 \< \xi \>^{2s} \, d\xi}.
\]
The sum over $m$ converges only when $s > 3/2$ in three dimensions.
This governs the dependence of the constant $C_s$ on $s$. Each $J_m$
gives one factor $\| f \|_s$, which completes the argument.
\end{proof}

Notice that the above result is sharp with respect to the range of $s$
for which it is valid. The Sobolev embedding from $H^s$ into continuous
functions is only valid when $s > 3/2$ in three dimensions, and it
does not in general make sense to sample a discontinuous function.

Let us note in passing that Lemma \ref{teo:sampling} generalizes
without difficulty to weighted norms.

\begin{lemma}\label{teo:sampling_weighted}
For all $s > 3/2$, $m \in \Z$, and $f \in H^s_m(\R^3)$,
\[
\| f - TSf \|_{L^2_m} \leq C_s \left( \frac{N}{L} \right)^{-s} \, \| f \|_{s,m}.
\]
\end{lemma}
\begin{proof}
  It is sufficient to exhibit a bandlimited multiplier $\xi_m(x)$ equivalent to $\< x \>^m$ in the sense that
\[
C_1 \< x \>^m \leq \xi_m(x) \leq C_2 \< x \>^m,
\]
for then weighted norms can be expressed with $\xi_m$ instead, and
\[
\xi_m (f - TSf) = (Id- TS) \xi_m f.
\]
Then the conclusion would follow from an application of Lemma
\ref{teo:sampling} to $\xi_m f$.

Such a function $\xi_m(x)$ can be constructed by convolving $\< x
\>^m$ by some appropriate nonnegative, band-limited kernel. For
example,
\[
\xi_m(x) = \< x \>^m * \prod_{j=1}^3
(\mbox{sinc}(\frac{x_j}{a}))^{2|m|+4}
\]
will do for $m \ne 0$, provided $a > 0$ is chosen so that the band limit of the
kernel is compatible with the sampling of $S$.
\end{proof}

We will frequently need to switch from discrete to continuous
norms, using the following result.

\begin{corollary}\label{teo:sampling2}
Let $f \in H^s(\R^3)$, for $s > 3/2$. Then
\[
\| f(x_j) \|_{\ell^2} \leq C_s \cdot \| f \|_{H^s}.
\]
\end{corollary}
\begin{proof}
  Decompose $TSf$ as $f + (TSf - f)$, and estimate in $L^2$. By Lemma~\ref{teo:sampling} and the sampling theorem,
\[
\| f(x_j) \|_{\ell^2} \leq \| f \|_{L^2} + C_s \cdot \left( \frac{N}{L} \right)^{-s} \| f \|_{H_s} \leq C \cdot \| f \|_{H^s}.
\]
\end{proof}

Let us collect some more background results. We have already hinted
at the fact that $(- \Delta)^{-1}$, defined as the convolution with
$G(x)$, is not bounded on $L^2$, or on any $H^s$ for that matter.
But it is bounded between \emph{weighted} Sobolev spaces. Let us
recall the following classical result from \cite{JenKat}.

\begin{lemma}\label{teo:JK} {\bf (Jensen-Kato)}
  The kernel $G(x-y) = \frac{1}{4 \pi |x-y|}$ maps boundedly
  $H^{-1}_m(\R^3)$ to $H^{1}_{-m'}(\R^3)$, and is in addition
  Hilbert-Schmidt from $L^2_m$ to $L^2_{-m'}$, provided
\[
m,m' > \frac{1}{2} \quad \text{and } \quad m+m' > 2.
\]
\end{lemma}


A similar property holds for the discretized kernel $\tilde{G}$.
Weigthed discrete $\ell^2$ spaces, relative to the grid $x_j$, are
defined as
\[
\ell^2_m = \{ u_j: \< x_j \>^{m} u_j \in \ell^2 \}.
\]

\begin{lemma}\label{teo:discJK}
  The kernel $\tilde{G}(x_j-x_k)$ defined in equation
  (\ref{eq:discfun}) maps boundedly $\ell^2_m$ to $\ell^2_{-m'}$
  provided
\[
m,m' > \frac{1}{2} \quad \text{and } \quad m+m' > 2.
\]
The operator norm of $\tilde{G}$ is uniform in $N/L \geq 1$.
\end{lemma}
\begin{proof}
The Hilbert-Schmidt norm of $\tilde{G}$ between $\ell^2_m$ and $\ell^2_{-m'}$ is
\[
\| \tilde{G} \|^2_{HS} = \sum_{j,k} |\tilde{G}(x_j-x_k)|^2 \< x_j
\>^{-2m'} \< x_{k} \>^{-2m}.
\]
The diagonal part $D$ of the sum, for $j = k$, is
\[
D = \frac{1}{(2 \pi)^2} \left( \frac{L}{N} \right)^4 \sum_{j} (1 + |j\frac{L}{N}|^2)^{-(m+m')}
\]
and can be compared to the integral
\[
D \leq C \cdot \frac{L}{N} \int_{\R^3} (1 + |x|^2)^{-(m+m')} \, dx,
\]
which converges when $m + m' > 3/2$. The off-diagonal part $OD$, for
$j \ne k$, is bounded by
\[
OD \leq \left( \frac{L}{N} \right)^6 \sum_{j,k: j \ne k} \frac{1}{|j\frac{L}{N} -
  k\frac{L}{N}|^2} (1 + |j\frac{L}{N}|^2)^{-m'} (1 + |k\frac{L}{N}|^2)^{-m}.
\]
(because $|Si(x)| < 2$.) Its continuous counterpart is
\[
OD \leq C \cdot \int_{\R^3} \int_{\R^3} \frac{1}{|x-y|^2} (1+|x|^2)^{-m}
(1+|y|^2)^{-m'} \, dx \, dy.
\]
The double integral is the Hilbert-Schmidt norm of $G$, squared,
between $L^2_m$ and $L^2_{-m'}$, and is bounded by lemma \ref{teo:JK}.
\end{proof}

We can now turn to the proof of the main result, theorem
\ref{teo:main}. The Birman-Schwinger operator is defined as
\[
K = U(x) (-\Delta^{-1}) U(x), \qquad \qquad U(x) = \phi(x)^\beta.
\]
The approximation $\tilde{K}$ is defined as
\[
\tilde{K} = T \tilde{U}(x_j) \tilde{G}_c \tilde{U}(x_j) S,
\]
where the subscript $c$ in $\tilde{G}_c$ indicates that the
convolution is not over the whole infinite grid $L/N \times \Z^3$, but
is a circular convolution over the finite cubic array $(-L/2 \, : \, L/N \, : \,
(L/2-L/N))^3$ (in Matlab notation.) Let us denote this bounded grid
by $\Box_N$.

\begin{proof}[Proof of theorem \ref{teo:main}.]

Let us divide the proof into four successive approximation steps, using the triangle inequality.

\begin{enumerate}
\item Let us first show that
\[
(I) = \| T \tilde{U}(x_j) (\tilde{G}_c - \tilde{G}) (\tilde{U}(x_j) S f) \|_{L^2}
\]
is adequately small in the sense of theorem \ref{teo:main}. By
Shannon's sampling theorem ($ST$ is the orthogonal projector onto
\Ran $\, S$,) we can rewrite
\[
(I) \leq \| \tilde{U}(x_j) (\tilde{G}_c - \tilde{G}) (\tilde{U}(x_j) f(x_j)) \|_{\ell^2}.
\]
We denote the operation of folding back a grid point by periodicity
onto the grid $\Box_N$ as follows:
\[
\lfloor x_j \rceil \equiv (x_j + (1,1,1)\frac{L}{2}) \mod L - (1,1,1)\frac{L}{2}.
\]
The discrepancy between the two types of convolution is
\[
(\tilde{G}_c - \tilde{G}) g(x_j) = \sum_{k \in \Z^3} g(x_k) \tilde{G}(x_j - x_k) - \sum_{x_k \in \Box_N} g(x_k) \tilde{G}( \lfloor x_j - x_k \rceil).
\]
Let us introduce the intermediate quantity
\[
\tilde{G}_b g(x_j) = \sum_{x_k \in \Box_N} g(x_k) \tilde{G}( x_j - x_k ),
\]
where the subscript $b$ stands for `bounded convolution'. The first contribution is
\begin{align*}
(I_A) &= \| \tilde{U}(x_j) (\tilde{G}_b - \tilde{G}) (\tilde{U}(x_j) f(x_j)) \|_{\ell^2}, \\
&= \| \tilde{U}(x_j) \sum_{x_k \notin \Box_N} \tilde{G}(x_j - x_k) \tilde{U}(x_k) f(x_k) \|_{\ell^2}, \\
&\leq \sup_j |\tilde{U}(x_j) \< x_j \>^2 | \cdot \| \tilde{G} \|_{\ell^2_{1} \to \ell^2_{-2}} \cdot \sup_{x_j \notin \Box_N} |\tilde{U}(x_j) \< x_j \>| \cdot \| f(x_j) \|_{\ell^2}.
\end{align*}
The first factor is obviously bounded by equation
(\ref{eq:assum_petv}), the second factor is bounded by Lemma
\ref{teo:discJK}, the third factor is less than $\epsilon + C \cdot
Le^{-\delta L}$ by equation (\ref{eq:assum_petv}), and the fourth
factor is less than $C \cdot \| f \|_{H^s}$ by Corollary
\ref{teo:sampling2}.

The second contribution is
\[
(I_B) = \| \tilde{U}(x_j) (\tilde{G}_c - \tilde{G}_b) (\tilde{U}(x_j) f(x_j)) \|_{\ell^2}.
\]
The kernels $\tilde{G}_c$ and $\tilde{G}_b$ differ only when
$x_j \notin \Box_N$, and in that case we have the estimate
\[
|\tilde{G}(x_j) - \tilde{G}( \lfloor x_j \rceil)| \leq \frac{1}{L} \chi_{x_j \notin \Box_N}(j,k).
\]
Therefore
\begin{align*}
(I_B) &\leq \frac{1}{L} \, \left\| \tilde{U}(x_j) \mathop{\sum_{x_k \in \Box_N}}_{x_k \notin x_j + \Box_N} \tilde{U}(x_k) f(x_k) \right\|_{\ell^2}, \\
&\leq \frac{1}{L} \left\| \, \tilde{U}(x_j) \cdot \sqrt{ \left( \frac{L}{N} \right)^3 \mathop{\sum_{x_k \in \Box_N}}_{x_k \notin x_j + \Box_N} |\tilde{U}(x_k)|^2} \; \right\|_{\ell^2} \cdot \| f(x_j) \|_{\ell^2}.
\end{align*}
The quantity underneath the square root can be bounded, up to a
multiplicative constant, by
\begin{align*}
\int_{x \notin B_{x_j}(L/2)} e^{-2 \delta |x|} \, dx &\leq \int_{x \notin B_{0}(L/2 - |x_j|)} e^{-2 \delta |x|} \, dx. \\
&\leq C \cdot \< \frac{L}{2} - |x_j| \>^2 e^{-2 \delta (\frac{L}{2} - |x_j|)}.
\end{align*}
With this bound, $(I_B)$ becomes
\begin{align*}
(I_B) &\leq \frac{C}{L} \sqrt{\left( \frac{L}{N} \right)^3 \sum_{x_j \in \Box_N} e^{-\delta |x_j}  \< \frac{L}{2} - |x_j| \> e^{-\delta (\frac{L}{2} - |x_j|)} }  \cdot \| f(x_j) \|_{\ell^2} \\
&\leq \frac{C}{L} \sqrt{\int_{B_0(L)} e^{-\delta|x|}  \< \frac{L}{2} - |x| \> e^{-\delta (\frac{L}{2} - |x|)} \, dx}  \cdot \| f(x_j) \|_{\ell^2}, \\
&\leq CL e^{-\delta L/4} \| f(x_j) \|_{\ell^2}.
\end{align*}
As before, the $\ell^2$ norm of $f(x_j)$ can be bounded, up to a
constant, by the $H^s$ norm of $f$ (Corollary \ref{teo:sampling2}).
This shows that $(I_A)$ and $(I_B)$ are both within the bounds of
equation (\ref{eq:main}).
\item Let us now study the difference
\[
(II) = \| T \tilde{U}(x_j) \tilde{G} (\tilde{U}(x_j) S f - S U(x) f) \|_{L^2}
\]
We have already observed that $T \tilde{U}(x_j)$ is bounded between
$\ell^2_{-2}$ and $L^2$. We also know from Lemma \ref{teo:discJK} that
$\tilde{G}$ is bounded from $\ell^2_1$ to $\ell^2_{-2}$. Hence, it
suffices to show that the following quantity is adequately small:
\[
\| (\tilde{U}(x_j) S - S U(x)) f \|^2_{\ell^2_{1}} = \left( \frac{L}{N} \right)^3 \sum_{j} \< x_j \>^2 |\tilde{U}(x_j) - U(x_j)|^2 |f(x_j)|^2.
\]
By equation (\ref{eq:assum_petv}) and Corollary \ref{teo:sampling2},
we can bound this expression by $C \epsilon^2 \| f \|_{H^s}^2$.
\item The third contribution
\[
(III) = \| (T \tilde{U}(x_j) - U(x) T) \tilde{G} S U(x) f \|_{L^2}
\]
can be bounded analogously. We know that $SU$ maps
$\ell^2$ to $\ell^2_{2}$ boundedly, and $\tilde{G}$ maps $\ell^2_2$ to
$\ell^2_{-1}$ boundedly. It remains to show that
\[
\| (T \tilde{U}(x_j) - U(x) T) g \|_{L^2}
\]
is small in proportion to $g \in \ell^2_{-1}$. By the sampling theorem, this is also
\[
\| (\tilde{U}(x_j) - S U(x) T) g(x_j) \|_{\ell^2} = \| (\tilde{U}(x_j) - U(x_j)) g(x_j) \|_{\ell^2}
\]
Again, equation (\ref{eq:assum_petv}) allows to bound this quantity by
$C \epsilon \| g(x_j) \|_{\ell^2_{-1}}$.
\item The last, and perhaps most important contribution, is
\[
(IV) = \| U (T \tilde{G} S + \Delta^{-1}) U f \|_{L^2}.
\]
Obviously, multiplication by $U(x)$ is bounded from $H^s$ to $H^s_2$,
as well as from $H^s_{-2}$ to $H^s$, for all $s \geq 0$. So it
suffices to show that
\[
\| (T \tilde{G} S + \Delta^{-1}) g \|_{L^2_{-2}} \leq C_s \cdot \left(
  \frac{N}{L} \right)^{-s} \cdot \| g \|_{H^s_2},
\]
for all $s > 3/2$. It is key to notice that
\[
T \tilde{G} S g = - \Delta^{-1} T S g,
\]
which follows from applying Lemma \ref{teo:exact_inverse_laplacian} to
$T S g$. Since by lemma \ref{teo:JK}, $\Delta^{-1}$ is bounded between
$L^{2}_2$ and $L^2_{-2}$, it is enough to check that
\[
\| TS g - g \|_{L^2_{2}} \leq C_s \cdot \left(
  \frac{N}{L} \right)^{-s} \cdot \| g \|_{H^s_2}.
\]
This is precisely the content of lemma \ref{teo:sampling_weighted}.
The proof is complete.
\end{enumerate}
\end{proof}

It is very likely that Theorem \ref{teo:main} could actually be
formulated with exponential decay in $N/L$, but we have been
unable to extend the argument. In order to do so, it would be
necessary to prove analyticity of the soliton $\phi(x)$, which we
believe is still an open problem in three dimensions. Our numerical
experiments strongly support that conjecture (see
Section~\ref{sec:results}).

\subsection{Accuracy of eigenvalues}

Let us remark that accuracy of the discretization depends crucially on
the smoothness of the function to which the operator is applied. This
follows not only from sampling issues, but also from the choice of
quadrature weigths $\tilde{G}(x_j)$ used to invert the laplacian. As a
consequence, only the eigenvalues of $K$ corresponding to very smooth
eigenfunctions will be computed to high accuracy.

Since $K$ is compact and self-adjoint on $L^2(\R^3)$, let us write
its spectral decomposition as
\[
K e_j = \lambda_j e_j, \qquad \< e_j, e_k \> = \delta_{jk}.
\]
Its discrete counterpart $\tilde{K}$ is also self-adjoint and compact
(of finite rank), but only on the space $B_{N/L}$. We have
\[
\tilde{K} \tilde{e}_j = \tilde{\lambda}_j \tilde{e}_j, \qquad \<
\tilde{e}_j, \tilde{e}_k \> = \delta_{jk}.
\]
Any $f \in B_{N/L}$ can therefore be expanded as $f = \sum_j \< f,
\tilde{e}_j \> \tilde{e}_j$.

The following result about accuracy of eigenvalues is mostly a
consequence of theorem \ref{teo:main}. Let $\epsilon_{L,N}$ denote the
small factor in the right-hand side of equation (\ref{eq:main}), namely
\[
\epsilon_{L,N} = C_{s} \cdot \left[ \epsilon + L e^{-\delta L/4} +
  \left( \frac{N}{L} \right)^{-s} \right].
\]

\begin{corollary}
If $\lambda_j$ is simple, then
\begin{equation}\label{eq:approxegv}
|\lambda_j - \tilde{\lambda}_j| \leq \frac{\epsilon_{L,N} \| e_j \|_{H^s}}{|\< e_j, \tilde{e}_j \>|},
\end{equation}
and, for $\epsilon_{L,N}$ sufficiently small,
\[
\| e_j - \tilde{e}_j \|_{L^2} \leq 2 \frac{\epsilon_{L,N} \| \tilde{e}_j \|_{H^s}}{d_j},
\]
where $d_j = \min_{k \ne j} |\lambda_j - \lambda_k|$.

If $\lambda_j$ has multiplicity $p > 1$, denote by $\Pi_j$ the
orthoprojector onto the $j$-th eigenspace, and $\tilde{\Pi}_j$ its numerical
counterpart. Then, for all $1 \leq m \leq p$,
\begin{equation}\label{eq:approxegv2}
|\lambda_j - \tilde{\lambda}_{j,m}| \leq \frac{\epsilon_{L,N}
  \sqrt{\sum_{n=1}^p \| e_{j,n} \|^2_{H^s}}}{\| \Pi_j \tilde{e}_{j,m}
  \|_{L^2}},
\end{equation}
and, for $\epsilon_{L,N}$ sufficiently small,
\[
\| \Pi_j - \tilde{\Pi}_j \|_{HS} \leq 2 \frac{\epsilon_{L,N} \sqrt{ \sum_{n=1}^p \| \tilde{e}_{j,n} \|^2_{H^s}}}{d_j}.
\]
The above norm is the Hilbert-Schmidt (HS) norm.
\end{corollary}
\begin{proof}
  The proof is loosely related to the argument behind Gershgorin's
  circle theorem.
\begin{itemize}
\item Take $\lambda_j$ simple. Let $\Pi = TS$ be the orthogonal projection from $L^2$ to
  $B_{N/L}$. We can then expand
\[
\Pi e_j = \sum_k \theta_{j,k} \tilde{e}_k, \qquad \theta_{j,k} = \< e_j, \tilde{e}_k \>.
\]
Consider the relation
\begin{equation}\label{eq:proofegv1}
\tilde{K} e_j = \lambda_j e_j + r_j,
\end{equation}
where by theorem \ref{teo:main}, the remainder $r_j = (\tilde{K} - K)
e_j$ obeys
\begin{equation}\label{eq:rjsmall}
\| r_j \|_{L^2} \leq \epsilon_{L,N} \| e_j \|_{H^s}.
\end{equation}
We can project equation (\ref{eq:proofegv1}) onto $B_{N/L}$ and expand
it on the basis $\tilde{e}_k$, which gives
\begin{equation}\label{eq:proofegv2}
\theta_{j,k} \tilde{\lambda}_k = \theta_{j,k} \lambda_j + \< r_j, \tilde{e}_k \>.
\end{equation}
For $k = j$, we can bound
\[
|\lambda_j - \tilde{\lambda}_j| \leq \frac{|\< r_j, \tilde{e}_j \>|}{|\theta_{j,j}|},
\]
which is exactly (\ref{eq:approxegv}) after using the estimate
(\ref{eq:rjsmall}) on $r_j$.

In order to obtain the estimate for the eigenfunctions, we should instead consider
\[
K \tilde{e}_j = \tilde{\lambda}_j \tilde{e}_j + \tilde{r}_j,
\]
with $\| \tilde{r}_j \|_{L^2} \leq \epsilon_{L,N} \| \tilde{e}_j
\|_{H^s}$, and this time expand it on the basis $e_k$,
\begin{equation}\label{eq:newthetajk}
\theta_{j,k} \lambda_k = \theta_{j,k} \tilde{\lambda}_{j} + \< \tilde{r}_j, e_k \>.
\end{equation}
For $k = j$, it would give the same estimate as (\ref{eq:approxegv}),
but with $\tilde{e}_j$ substituted for $e_j$. For $k \ne j$, we can
rewrite equation (\ref{eq:newthetajk}) as
\begin{equation}\label{eq:thetajk}
|\theta_{j,k}| \leq \frac{|\< \tilde{r}_j, e_k \>|}{|\tilde{\lambda}_j - \lambda_k|}.
\end{equation}
We are not quite finished since the eigevalue gap, at the denominator,
is measured with $\tilde{\lambda_j}$ instead of $\lambda_j$. Put
$\delta_j = |\lambda_j - \tilde{\lambda}_j|$. Let us also introduce
$\tilde{d}_j = \min_{k \ne j} |\tilde{\lambda}_j - \lambda_k|$, and
observe that $\tilde{d}_j \geq d_j - \delta_j$. Then, squaring
(\ref{eq:thetajk}) and summing over $k \ne j$, we get
\[
\sum_{k \ne j} |\theta_{j,k}|^2 \leq \frac{\| \tilde{r}_j \|^2_{L^2}}{\tilde{d}^2_j},
\]
and, as a result,
\begin{equation}\label{eq:thetajj}
|\theta_{j,j}|^2 \geq 1 - \frac{\| \tilde{r}_j \|^2_{L^2}}{\tilde{d}^2_j}.
\end{equation}
Our estimate for $\delta_j$, equation (\ref{eq:approxegv}), can
therefore be improved to
\[
\delta^2_j \leq \frac{\epsilon^2_{L,N} \| e_j \|^2_{H^s}}{1 -
  \frac{\epsilon^2_{L,N} \| \tilde{e}_j \|^2_{H^s}}{(d_j -
    \delta_j)^2}}.
\]
It is not hard to see that, as $\epsilon_{L,N}$ gets small,
$\delta_j$, defined here implicitly, decreases to zero. Take $L$, $N$
and $1/\epsilon$ sufficiently large so that $\delta_j \leq d_j/4$.
Going back to equation (\ref{eq:thetajj}), we obtain
\[
|\theta_{j,j}|^2 \geq 1 - 2 \frac{\| \tilde{r}_j \|^2_{L^2}}{d^2_j}.
\]
The estimate for eigenvectors follows since
\[
\| e_j - \tilde{e}_j \|^2 = 2 - 2 \theta_{jj} \leq 2 - 2 \theta_{jj}^2.
\]
\item Assume now that $\lambda_j$ has multiplicity $p > 1$. The previous
  argument goes through, with the following modifications. The change of
  basis coefficients are now $\theta_{jn,k} = \< e_{jn}, \tilde{e}_k
  \>$, where $n = 1, \ldots, p$.  Equation (\ref{eq:proofegv2}) becomes
\begin{equation}\label{eq:proofegv3}
\theta_{jn,k} \tilde{\lambda}_k = \theta_{jn,k} \lambda_j + \< r_{j,n}, \tilde{e}_k \>.
\end{equation}
Let us study the discrepancy between $\lambda_j$ and $\tilde{\lambda}_{jm}$,
for some $m$. We need to take a linear combination of equation
(\ref{eq:proofegv3}) with the weights
\[
\alpha_n = \frac{\theta_{jn,jm}}{\sqrt{\sum_n |\theta_{jn,jm}|^2}}.
\]
We then obtain the inequality
\[
|\lambda_j - \tilde{\lambda}_{jm}| \leq \frac{\sum_n \alpha_n |\<
  r_{j,n}, \tilde{e}_{jm} \>|}{\sqrt{\sum_n |\theta_{jn,jm}|^2}}.
\]
Observe that $\| \Pi_j \tilde{e}_{jm} \|_{L^2} = \sqrt{\sum_n
  |\theta_{jn,jm}|^2}$, and the numerator can be bounded by
Cauchy-Schwarz (in $n$) to yield (\ref{eq:approxegv2}).

As for the eigenfunctions, Hilbert-Schmidt norms and trace inner
products are to be substituted for their $L^2$ scalar counterpart. More precisely,
\begin{align*}
\| \Pi_j - \tilde{\Pi}_j \|^2_{HS} &= \mbox{Tr} (\Pi_j - \tilde{\Pi}_j )^2, \\
&= \mbox{Tr} \Pi_j + \mbox{Tr} \tilde{\Pi}_j - 2 \mbox{Tr} \Pi_j \tilde{\Pi}_j, \\
&= 2p - 2 \mbox{Tr} \Pi_j \tilde{\Pi}_j,
\end{align*}
so we need to find a bound on the latter quantity. Observe first that
\[
\mbox{Tr} \Pi_j \tilde{\Pi}_j = \sum_n \sum_{m=1}^p |\theta_{jn,jm}|^2.
\]
The change of basis coefficients are normalized in the sense that, for
all $n = 1, \ldots, p$,
\[
1 = \sum_k |\theta_{jn,k}|^2 = \sum_m |\theta_{jn,jm}|^2 + \sum_{k \ne j } |\theta_{jn,k}|^2.
\]
Summing over $n$, we get
\begin{equation}\label{eq:proofegv4}
p = \mbox{Tr} \Pi_j \tilde{\Pi}_j + \sum_n \sum_{k \ne j } |\theta_{jn,k}|^2.
\end{equation}
The bound on $\theta_{jn,k}$ for $k \ne j$ can be obtained as above, and is
\begin{equation}\label{eq:proofegv5}
|\theta_{jn,k}| \leq \frac{|\< \tilde{r}_{jn}, e_k
  \>|}{|\tilde{\lambda}_{jn} - \lambda_k|}.
\end{equation}
It is natural to define $\tilde{d}_j = \min_{k,n}
|\tilde{\lambda}_{jn} - \lambda_k|$. The way to bound $\tilde{d}_j$
from below by, say, $3 d_j/4$ is very analogous to the simple case,
and based on some adequate control of the size of the denominator in
equation (\ref{eq:approxegv2}). This can be obtained from
\[
  1 = \| \Pi_j \tilde{e}_{jm} \|_2^2 + \sum_{k \ne jn} |\theta_{k,jm}|^2,
\]
and
\[
\sum_{k \ne jn} |\theta_{k,jm}|^2 \leq \sum_{k \ne jn} \frac{|\< r_k,
  \tilde{e}_{jm} \>|^2}{|\lambda_k - \tilde{\lambda}_{jm}|^2} =
\sum_{k \ne jn} \frac{|\< e_k, \tilde{r}_{jm} \>|^2}{|\lambda_k -
  \tilde{\lambda}_{jm}|^2} \leq \frac{\| \tilde{r}_{jm}
  \|_{L^2}^2}{\tilde{d}^2_j}.
\]
(we have used $K = K^*$ and $\tilde{K} = \tilde{K}^*$ for
band-limited functions.) Eventually, we use theorem (\ref{teo:main})
one more time to bound
\begin{equation}\label{eq:proofegv6}
\| \tilde{r}_{jn} \|_{L^2} \leq \epsilon_{L,N} \| \tilde{e}_{jn} \|_{H^s}.
\end{equation}
The desired estimate is obtained by combining the intermediate
inequalities (\ref{eq:proofegv4}) to (\ref{eq:proofegv6}). This
concludes the proof.
\end{itemize}
\end{proof}

Note that, in the formulation of the corollary or in its
justification, it is nowhere necessary to obtain bounds on errors in
computing \emph{other} eigenvalues $\lambda_k \ne \lambda_j$. This is
the sense in which our argument is close to Gershgorin's theorem: we
only need one row of the matrix $\< e_j, \tilde{K} e_k \>$, namely,
the $j$-th row.

\section{Numerical results and discussion}\label{sec:results}

\begin{figure}[htb]
\includegraphics[width=8cm]{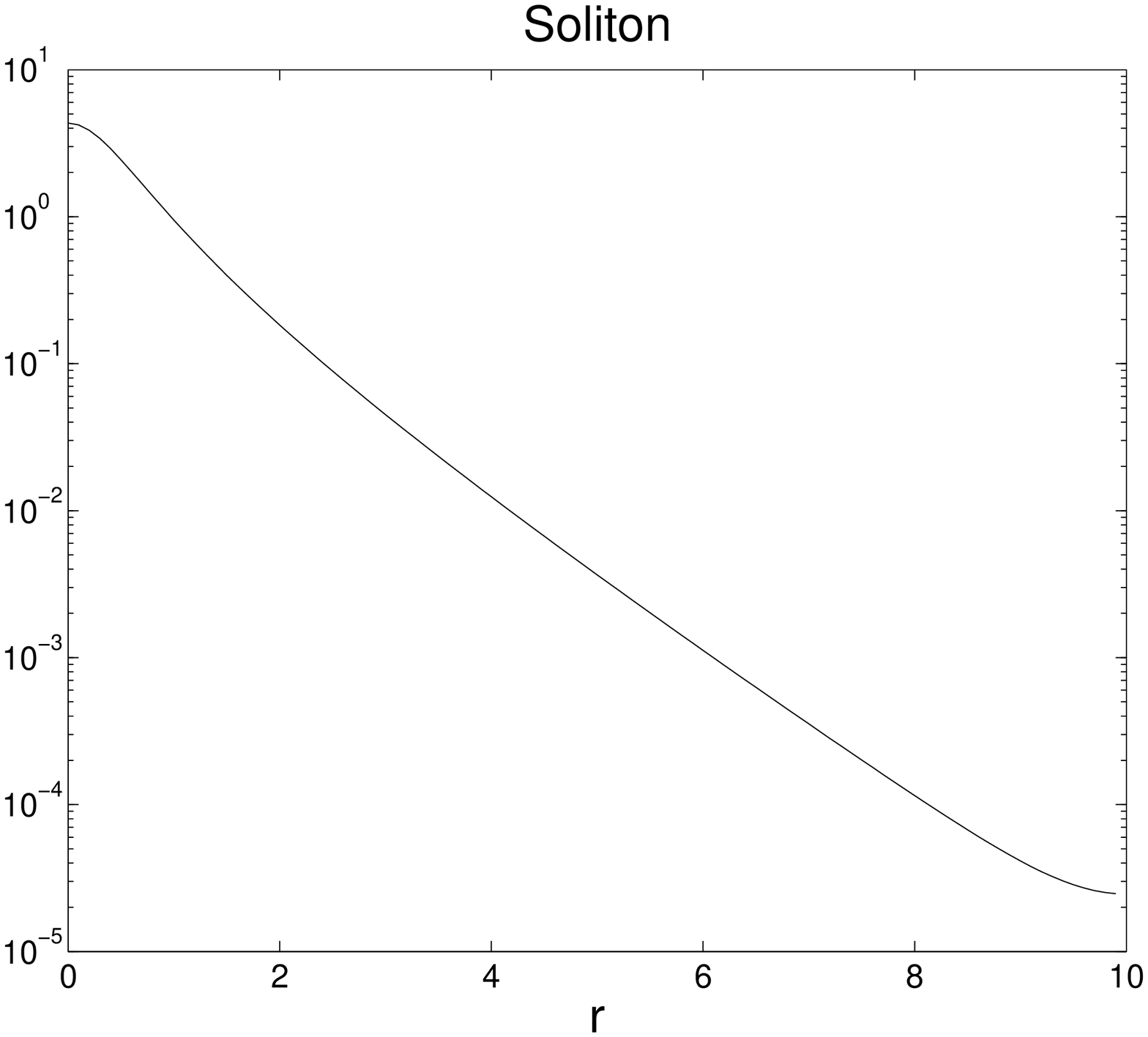}\quad
\includegraphics[width=8cm]{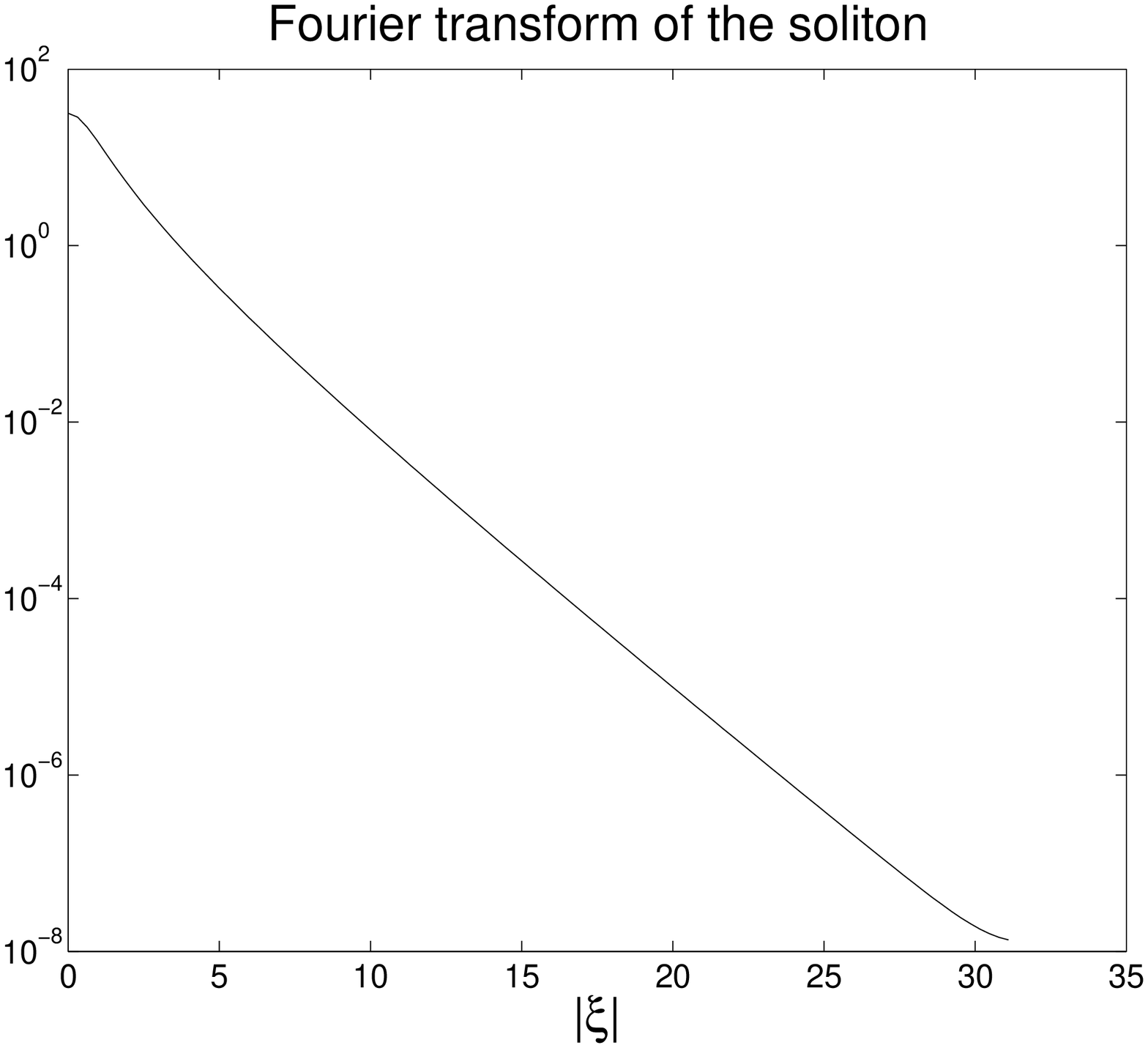}
\caption{Left: the soliton $\phi(r)$ in space. Right: the soliton $\hat\phi(|\xi|)$ in frequency. Both depend only on the radial coordinate.}
  \label{fig:soliton}
\end{figure}

\begin{figure}[htb]
\begin{center}\includegraphics[width=12cm]{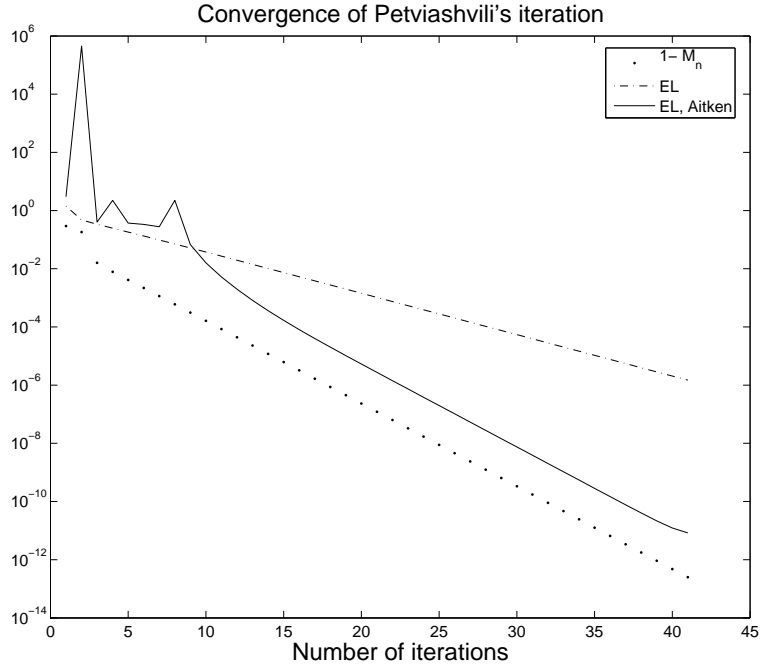}\end{center}
\caption{Convergence of Petviashvili's iteration. The x-axis is the number of iteration. Dotted line: $1 - M_n$, where $M_n$ is the Petviashvili constant. Dash-dotted line: Euler-Lagrange remainder, $ \| - \Delta
  \phi_n + \phi_n - \phi_n^{2 \beta + 1} \|/\| \phi_n \|$, with norms
  in $L^2$. Solid line: Same remainder, but for the Aitken iterate
  $\phi_n^A$.}\label{fig:petv_conv}
\end{figure}

Our numerical approximation to the soliton $\phi(x)$, obtained by
Petviashvili's iteration, is plotted in log-scale in Figure
\ref{fig:soliton}. For this experiment, we have taken $\beta = 1$, a
cube of sidelength $L = 20$ and a grid of size $N = 200$. Notice the
apparent exponential decay both in space and frequency. Figure
\ref{fig:petv_conv} checks the convergence of Petviashvili's
iteration, namely that $M_n \to 1$ and that $- \Delta \phi_n + \phi_n
- \phi_n^{2 \beta + 1} \to 0$ as $n \to \infty$. The constants
$R_{n,j}$ are not plotted; in all our tests they are at machine
accuracy, or below. In other words, the soliton is as radial as
allowed by the grid.

In the above setting we observe nice convergence of the error, as
measured by different indicators, to machine accuracy. This is not
always the case. For values of $L$ too big in comparison to $N$
(typically, for $L$ of the order of $N/4$ and above,) the
Petviashvili iteration possibly reaches acceptable error levels, but
then diverges away from the fixed point. A careful inspection of the
remainder $- \Delta \phi_n + \phi_n - |\phi_n^{2 \beta}|\phi_n$
indicates that this might be due the fifth eigenvalue $\lambda_5$ of
the linearized iteration operator $A$ becoming negative, in the
notations of Section~\ref{sec:petv_conv}.

Next, we show a plot of the two, resp. five largest eigenvalues of the
Birman-Schwinger operator $K_-$, resp. $K_+$, as a function of $\beta$
in the range $[\frac{2}{3},1]$. As mentioned earlier, $\lambda_2(K_-)$
is less than 1 for all values of $\beta$, but there exists a number
$\beta_*$ below which $\lambda_5(K_+) > 1$. This is the signature of
at least one eigenvalue in the gap of $L_+$, for $2/3 \leq \beta <
\beta_*$ (inspection of $\lambda_6$ reveals that there is only one
eigenvalue in the gap.) For this experiment, we used $L = 15$ and $N =
60$. Note that, in both cases, $\lambda_2 = \lambda_3 = \lambda_4$ is
triple, and $\lambda_5$ is simple, so we only show at most 3 curves.
Our numerical implementation correctly picks up the multiplicity,
exactly (to all 16 digits), and in all the cases that we have tried.

\begin{figure}[htb]
\includegraphics[width=8cm]{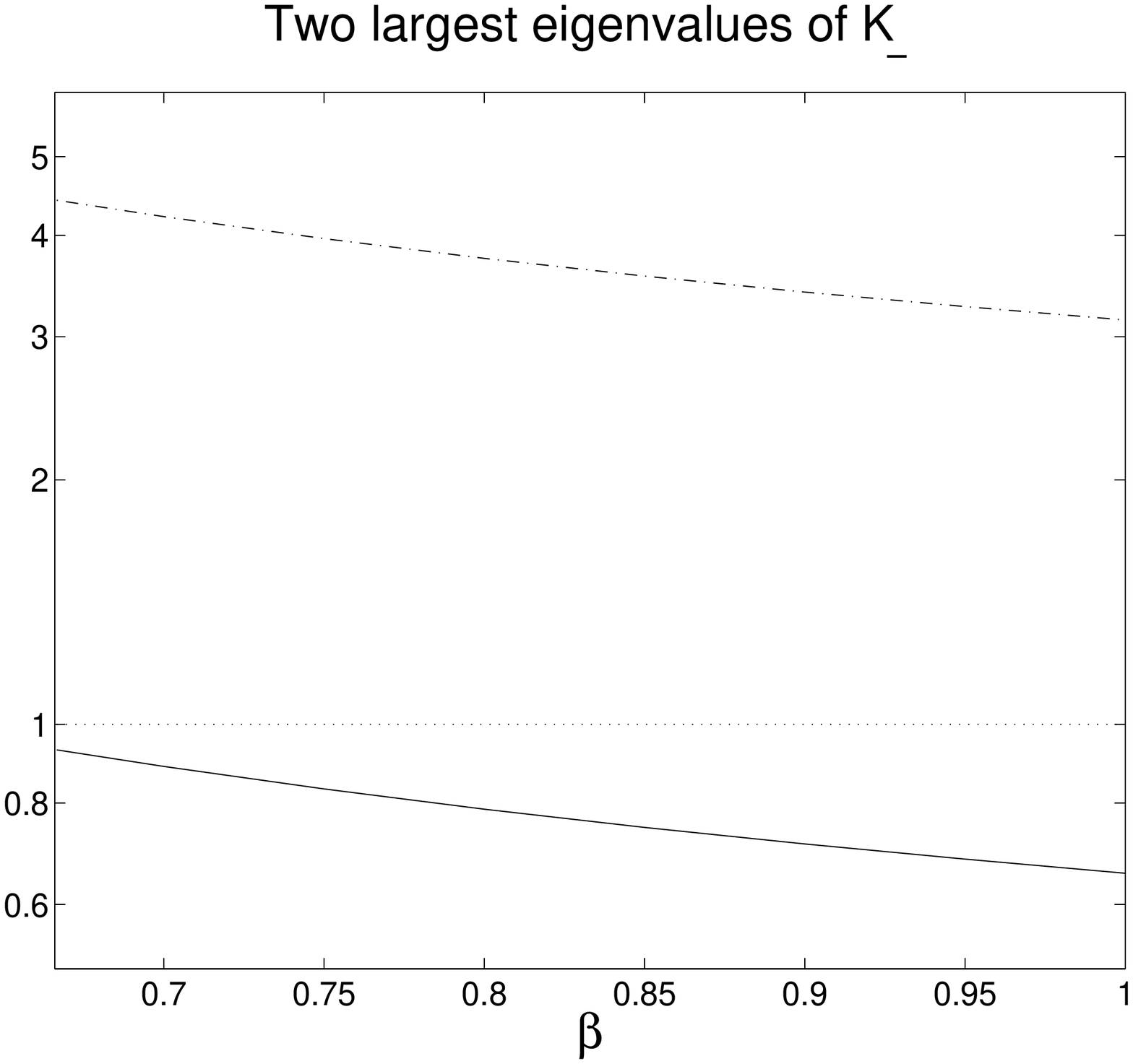}\quad
\includegraphics[width=8cm]{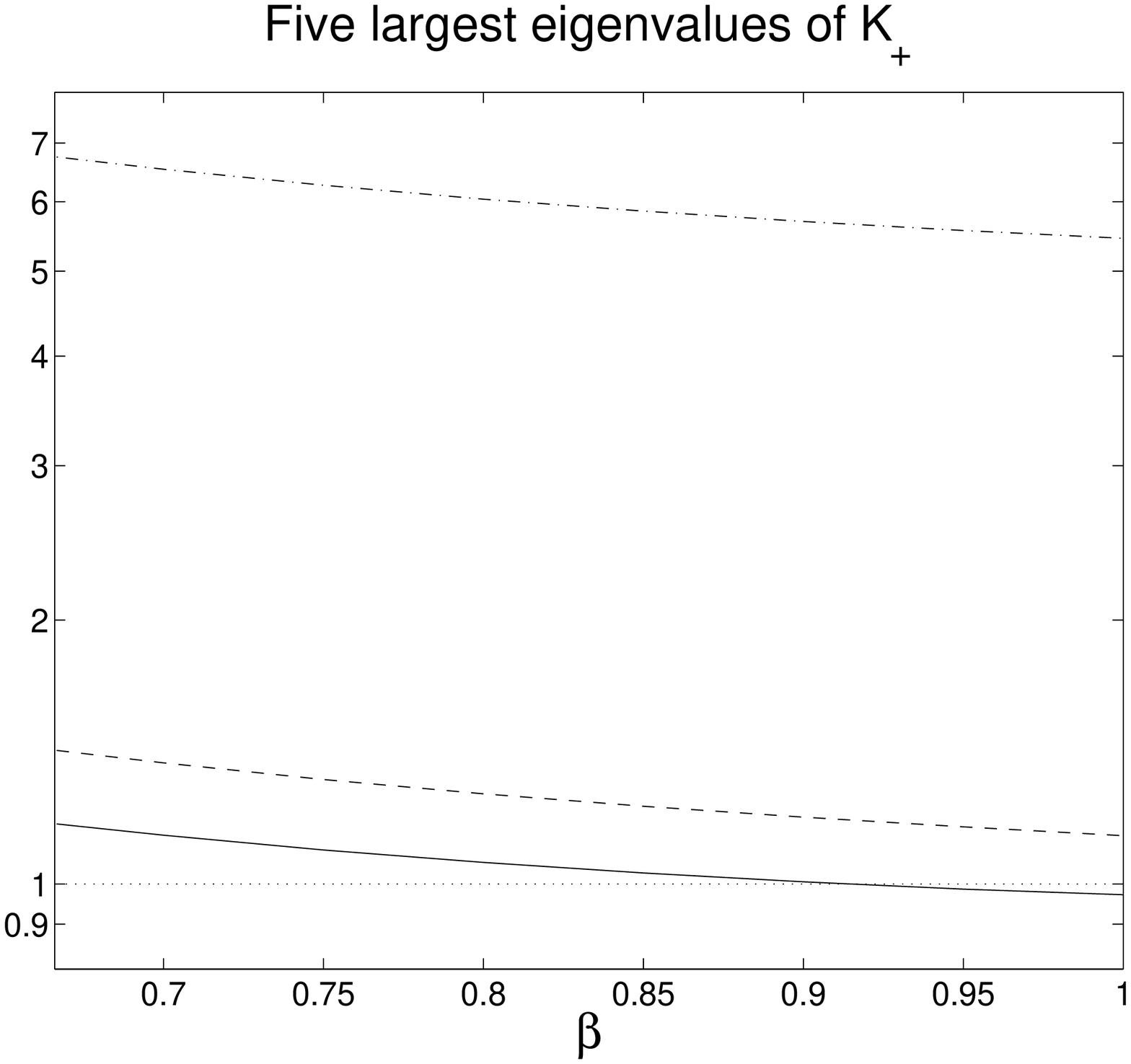}
\caption{Left: Two largest eigenvalues of $K_-$, as a function of $\beta$ in the range $[\frac{2}{3},1]$.  Right: Five largest eigenvalues of $K_+$, as a function of $\beta$. In both cases, the largest eigenvalue is simple, the second to fourth eigenvalues are one triplet, and the fifth eigenvalue is simple. The value one is indicated by the dotted line.  }
  \label{fig:lambdas}
\end{figure}

An accurate computation of the numerical value of the exceptional
exponent $\beta_*$ requires higher values of $L$ and $N$. On a 2005
standard desktop, we have tried $L = 25$ and $N = 200$. We then obtain
$\beta_*$ by interpolation of $\lambda_5(K_+)$ for different closeby
values of $\beta$, see table \ref{tbl:interp}. The confidence on
$\beta_*$, which we estimate to be about 8 digits, is directly related
to the level of accuracy of $\lambda_5(K_+)$. The latter is determined
by inspection of convergence as $L$ and $N$ increase. So the bounds
given in the main claim, in the introduction, are \emph{not} rigorous,
but merely serve as an indication of how accurate we believe our
algorithm is.

\begin{table}[ht!]
  \begin{center}
    \begin{tabular}{|cc|}
      \hline
      $\beta$ & $\lambda_5(K_+)$ \\
      \hline
      0.91395850 & 1.00000016477 \\
      0.91395875 & 1.00000006304 \\
      0.91395900 & 0.99999996130 \\
      0.91395925 & 0.99999985957 \\
      \hline
    \end{tabular}
  \end{center}
  \caption{Fifth eigenvalue of $K_+$, as a function of $\beta$ near $\beta_*$. Here $L = 25$ and $N = 200$. Each value took about one day to obtain. Cubic interpolation reveals $\beta_* \simeq 0.913958905 \pm 1e-8$.}
  \label{tbl:interp}
\end{table}

At this point the reader might wonder if, instead of a full
three-dimensional simulation, there exists a computational strategy
involving only the radial coordinate to compute both the soliton and
the eigenvalues of the Birman-Schwinger operators. We believe the
answer is positive, but will involve significantly different ideas
from the ones presented in this paper. In particular, spectral
accuracy will be more difficult to obtain. We leave the problem of
determining more digits of the constant $\beta_*$, likely through the
use of a one-dimensional method, as a challenge to the interested
reader.

Finally, the Matlab code we used to generate the figures and compute
an estimation of $\beta_*$ can be freely downloaded from
\begin{verbatim}http://www.acm.caltech.edu/~demanet/NLS/.\end{verbatim}

\bibliographystyle{amsplain}

\begin{thebibliography}{99}



\bibitem[BerCaz]{BerCaz} Berestycki, H., Cazenave, T.
{\em Instabilit\'e des \'etats stationnaires dans les \'equations de Schr\"odinger et de Klein-Gordon non lin\'eaires.}
C.\ R.\ Acad.\ Sci.\ Paris S\'er.~I Math. 293 (1981), no. 9, 489--492.

\bibitem[BerLio]{BerLio} H. Berestycki, P.L. Lions, {\em Nonlinear scalar field equations. I. Existence of a ground
state.} Arch.\ Rational Mech.\ Anal.\ 82 (1983), no.~4, 313--345.

\bibitem[BerLioPel]{BerLioPel} H. Berestycki, P.L. Lions, L.A.
  Peletier, {\em An ODE approach to the existence of positive
    solutions for semilinear problems in $\R^n$.} Indiana U. Math. J.
  30 (1981), no.~1, 141--157.


\bibitem[BusPer1]{BP1} Buslaev, V.\ S., Perelman, G.\ S. {\em Scattering for the
nonlinear Schr\"odinger equation: states that are close to a soliton.} (Russian)
Algebra i Analiz  4  (1992),  no.~6, 63--102;  translation in  St.\ Petersburg Math.\ J.~4
(1993),  no.~6, 1111--1142.


\bibitem[CazLio]{CazLio} Cazenave, T., Lions, P.-L. {\em Orbital stability
of standing waves for some nonlinear Schr\"odinger equations.} Comm.\ Math.\ Phys.~85 (1982), 549--561



\bibitem[Cuc]{Cuc} Cuccagna, S.  {\em Stabilization of solutions to nonlinear Schr\"odinger equations.}
Comm.\ Pure Appl.\ Math.~54  (2001),  no.~9, 1110--1145.

\bibitem[Eigs]{eigs}
\begin{verbatim}http://www.mathworks.com/access/helpdesk/help/techdoc/ref/eigs.html,\end{verbatim}
and references therein.

\bibitem[ErdSch]{ErdSch} Erdo\smash{\u{g}}an, M. B., Schlag, W.
{\em Dispersive estimates for Schr\"{o}dinger operators in the
presence of a resonance and/or an eigenvalue at zero energy in
dimension three: II}, preprint 2005.

\bibitem[Flu]{Flug}  Fl\"ugge, S. {\em Practical quantum mechanics.}
Reprinting in one volume of Vols. I, II. Springer-Verlag, New
York-Heidelberg, 1974.







\bibitem[Gri]{Grill} Grillakis, M. {\em
Analysis of the linearization around a critical point of an infinite dimensional
Hamiltonian system.} Comm.\ Pure Appl.\ Math.~41 (1988), no.~6, 747--774.

\bibitem[GriShaStr1]{GSS1} Grillakis, M., Shatah, J., Strauss, W.
{\em Stability theory of
solitary waves in the presence of symmetry. I.}
J.\ Funct.\ Anal.~74  (1987),  no.~1, 160--197.

\bibitem[GriShaStr2]{GSS2} Grillakis, M., Shatah, J., Strauss, W.
{\em Stability theory of
solitary waves in the presence of symmetry. II.}
J.\ Funct.\ Anal.~94  (1990), 308--348.


\bibitem[HutPym]{HutPym} Hutson, V., Pym, J.S., {\em Applications of
  Functional Analysis and Operator Theory}, Academic Press, London, 1980.

\bibitem[JenKat]{JenKat} Jensen, A., Kato, T., {\em Spectral
    properties of Schr\"{o}dinger operators and time-decay of the wave
    functions,} Duke Math. J. 46 (1979), no.3,  583--611



\bibitem[Kwo]{Kwo} Kwong, M.\ K. {\em Uniqueness of positive solutions
of $\Laplace u - u +u^{p}=0$ in $\R^{n}$.} Arch.\ Rat.\ Mech.\ Anal.~65
(1989), 243--266.



\bibitem[PelSte]{PelSte} D. Pelinovsky, Y. Stepanyants, {\em
    Convergence of Petviashvili's iteration method for numerical
    approximation of stationary solutions to nonlinear wave
    equations.} SIAM J. Num. Anal. 42 (2004), no.~3, 1110--1127
%

\bibitem[Per2]{Pe2} Perelman, G. {\em On the formation of singularities in solutions
of the critical nonlinear Schrödinger equation.} Ann.\ Henri
Poincar\'e 2 (2001), no.~4, 605--673.






\bibitem[RodSchSof1]{RSS1} Rodnianski, I., Schlag, W., Soffer, A.
{\em Dispersive Analysis of Charge Transfer Models}, preprint 2002, to appear in CPAM

\bibitem[RodSchSof2]{RSS2} Rodnianski, I., Schlag, W., Soffer, A.
{\em Asymptotic stability of $N$-soliton states of NLS}, preprint 2003, submitted to CPAM


\bibitem[Sch]{schlag} Schlag, W. {\em Stable manifolds for orbitally unstable NLS.} preprint, 2004.







\bibitem[Str]{Str} Strauss, W.\ A. {\em Nonlinear wave equations.} CBMS Regional
Conference Series in Mathematics, 73. American Mathematical
Society, Providence, RI, 1989


\bibitem[SulSul]{SulSul} Sulem, C., Sulem, P.-L.
{\em The nonlinear Schr\"odinger equation. Self-focusing and wave collapse.}
 Applied Mathematical Sciences, 139. Springer-Verlag, New York, 1999.



\bibitem[Wei1]{Wei1} Weinstein, Michael I.
{\em Modulational stability of ground states of nonlinear Schr\"odinger equations.}
  SIAM J.\ Math.\ Anal.~16  (1985),  no.~3, 472--491.

\bibitem[Wei2]{Wei2} Weinstein, Michael I.
{\em Lyapunov stability of ground states of nonlinear dispersive evolution equations.}
Comm.\ Pure Appl.\ Math.~39  (1986),  no.~1, 51--67.

\end{thebibliography}

\noindent

\medskip\noindent
\textsc{Department of Applied and Computational Mathematics, 217-50 Caltech, \\  Pasadena, CA 91125, U.S.A.}\\
{\em email: }\textsf{demanet@acm.caltech.edu},

\medskip\noindent
\textsc{Department of Mathematics,
253-37 Caltech,\\ Pasadena, CA 91125, U.S.A.}\\
{\em email: }\textsf{schlag@its.caltech.edu}

\end{document}